\documentclass[11pt]{article}

\usepackage{color}
\usepackage[T1]{fontenc}
\usepackage[utf8]{inputenc} 

\newtheorem{myproposition}{Proposition}[section]
\newtheorem{mytheorem}[myproposition]{Theorem}
\newtheorem{mylemma}[myproposition]{Lemma}

\newtheorem{myconjecture}[myproposition]{Conjecture}
\newtheorem{mycorollary}[myproposition]{Corollary}

\newtheorem{myobservation}[myproposition]{Observation}

\newtheorem{myproblem}[myproposition]{Problem}
\usepackage{amsfonts}
\usepackage{amsmath}
\def\gr{\mathcal{G}}
\def\zet{\mathbb{Z}}
\def\gchi{\chi^\Sigma_g}
\newcommand{\qed}{\hfill \rule{.1in}{.1in}}

\makeatletter
\def\imod#1{\allowbreak\mkern10mu({\operator@font mod}\,\,#1)}
\makeatother

\begin{document}

% A short title is not required, but if needed use:
% \title[short title]{full title}
\title{Linear bounds on nowhere-zero group irregularity strength and nowhere-zero group sum chromatic number of graphs}

%\author{Marcin Anholcer\\
%{\small  Faculty of Informatics and Electronic Economy}\\
%{\small Pozna\'n University of Economics and Business}\\
%{\small Al. Niepodleg{\l}o\'sci 10, 61-875 Pozna\'n, Poland}\\
%}

\author{Marcin Anholcer$^1$, Sylwia Cichacz$^{2,3}$\footnote{This work was partially supported by the Faculty of Applied Mathematics AGH UST statutory tasks within subsidy of Ministry of Science and Higher Education.}, Jakub Przyby{\l}o$^2$\\
$^1$Pozna\'n University of Economics and Business\\
$^2$AGH University of Science and Technology,\\ 
al. A. Mickiewicza 30, 30-059 Krakow, Poland\\ %Krak\'ow, Poland\\
$^3$Faculty of Mathematics, Natural Sciences and Information Technologies,\\
 University of Primorska, Glagolja\v{s}ka 8, SI 6000 Koper, Slovenia}

%\email{m.anholcer@ue.poznan.pl}

%{\small Faculty of Applied Mathematics}\\
%{\small AGH University of Science and Technology}\\
%{\small Al. Mickiewicza 30, 30-059 Krak\'ow, Poland}\\

\maketitle

% For each additional author, add another set of
% \author, \address, and \email commands
%
%\author{}
%\address{}
%\email{}

%\subjclass[2000]{05C15,05C78}
%\keywords{Irregularity strength, Graph weighting, Graph labeling, Abelian group}

% \thanks entries are to acknowledge grants. You may combine
% all acknowledgments into one \thanks entry, or may use
% multiple \thanks entries. They generate footnotes without
% tags, so you must be explicit about which authors are
% thanking whom.
%\thanks{}

\begin{abstract}
We investigate the \textit{group irregularity strength}, $s_g(G)$, of a graph, i.e. the least integer $k$ such that taking any Abelian group $\gr$ of order $k$, there exists a function $f:E(G)\rightarrow \gr$ so that the sums of edge labels incident with every vertex are distinct. So far the best upper bound on $s_g(G)$ for a general graph $G$ was exponential in $n-c$, where $n$ is the order of $G$ and $c$ denotes the number of its components. In this note we prove that $s_g(G)$ %the upper bound 
is linear in $n$, namely not greater than $2n$. In fact, we prove a stronger result, as we additionally forbid the identity element of a group to be an edge label or the sum of labels around a vertex.
%do not use the identity element of groups neither on edges nor on vertices. 
We consider also locally irregular labelings where we require only sums of adjacent vertices to be distinct. For the corresponding graph invariant %${\gchi}^\star$ 
we prove the general upper bound: $\Delta(G)+{\rm col}(G)-1$ (where ${\rm col}(G)$ is the coloring number of $G$) in the case when we do not use the identity element as an edge label, and a slightly worse one if we additionally forbid it as the sum of labels around a vertex. 
In the both cases we also provide a sharp upper bound for trees and a constant upper bound for the family of planar graphs.
%This time we show that the upper bound is linear in $\Delta(G)$ and ${\rm col}(G)$ - KONKRETNIE. We also %provide a sharp upper bound in the case of trees and a constant upper bound for the family of planar graphs.
\end{abstract}

% Insert paper text here
\section{Introduction}

It  is a well known fact that in any simple graph $G$
there are at least two vertices of the same degree. The situation changes if we consider an edge labeling $f:E(G)\rightarrow \{1,\ldots,k\}$ and calculate so-called \textit{weighted} degree (or \textit{weight}) of each vertex $v$ as the sum of labels of all the edges incident to $v$. The labeling $f$ is called \textit{irregular} if the weighted degrees of all the vertices are distinct. The smallest value of $k$ that allows some irregular labeling is called the \textit{irregularity strength of $G$} and denoted by $s(G)$.

The problem of finding $s(G)$ was introduced by Chartrand et al. in \cite{ref_ChaJacLehOelRuiSab1} and investigated by numerous authors  \cite{ref_AigTri,ref_AmaTog,ref_FerGouKarPfe,ref_Leh,ref_Tog1}.  An upper bound $s(G)\leq n-1$ was proved for all graphs containing no isolated edges and at most
one isolated vertex, except for the graph $K_3$ \cite{ref_AigTri,ref_Nie}, where $n$ is the order of $G$. This is %a 
tight, %upper bound,
as exemplified e.g.  by the family of stars. % A better  upper bound can
It can however be improved  for graphs with sufficiently large   minimum  degree $\delta$.  The best published general result due to Kalkowski et al. (see \cite{ref_KalKarPfe1}) is $s(G)\leq 6n/\delta$. It was recently improved by Majerski and Przyby{\l}o (\cite{ref_Prz3}) for relatively dense graphs with sufficiently large minimum degree compared to $n$ ($s(G)\leq (4+o(1))n/\delta+4$ in this case). 

Jones combined the concepts of graceful labeling and modular edge coloring into a labeling called  a \emph{modular edge-graceful labeling} (\cite{ref_Jon,ref_JonKolOkaZha2, ref_JonZha}). He defined the \textit{modular edge-gracefulness} of a graph $G$ as the smallest integer $k(G)=k \geq n$ for which there exists an edge labeling $f:E(G)\rightarrow \zet_k$ such that the induced vertex labeling $f^\prime : V (G)\rightarrow \zet_k$ defined by
$$
f^\prime(v) =\sum_{u\in N(v)}f(uv)\mod k
$$
is one-to-one.

Assume $\gr$ is an Abelian group  of order $m\geq n$ with the operation denoted by $+$ and identity element $0$. For convenience we will write $ka$ to denote $a+a+\ldots+a$ (where element $a$ appears $k$ times), $-a$ to denote the inverse of $a$ and we will use $a-b$ instead of $a+(-b)$. Moreover, the notation $\sum_{a\in S}{a}$ will be used as a short form for $a_1+a_2+a_3+\ldots$, where $a_1, a_2, a_3, \ldots$ are all the elements of the set $S$. 
%Recall that any group element $\iota\in\gr$ of order 2 (i.e., $\iota\neq 0$ such that $2\iota=0$) is called an
%\emph{involution}. 

%The order of an element $a\neq 0$ is the smallest $r$ such that $ra=0$. It is well-known by Lagrange Theorem that $r$ divides $|\gr|$ \cite{ref_Gal}. Therefore every  group of odd order has no involution.

We consider an edge labeling $f:E(G)\rightarrow \gr$ inducing %leading us to the weighted degrees 
the \textit{weighted degree} of every vertex $v$ of $G$
defined as the sum (in $\gr$):
$$
w_f(v)=\sum_{u\in N(v)}f(uv).
$$
This we shall also call the \textit{weight} of $v$ or the \textit{sum at} $v$ and denote simply by $w(v)$ if this causes no ambiguities. 
The labeling $f$ of $G$ is called $\gr$-\textit{irregular}, if the resulting weighted degrees of all the vertices are pairwise distinct. The least positive integer $k$ such that for every Abelian group $\gr$ of order $k$ there exists a $\gr$-irregular labeling of $G$ is called the the \textit{group irregularity strength} of $G$ and denoted $s_g(G)$.

%The concept of $\gr$-irregular labeling is a generalization of modular edge-graceful labeling. In both cases the labeling $f$ is called \textit{$\gr$-irregular} if all the weighted degrees are distinct. However, the \textit{group irregularity strength} of $G$, denoted $s_g(G)$, is the smallest integer $s$ such that for every Abelian group $\gr$ of order $s$ there exists $\gr$-irregular labeling $f$ of $G$. Thus the following observation is true.
%
%\begin{myobservation}[\cite{ref_AnhCic2}]
%For every graph $G$ with no component of order less than $3$, $k(G)\leq s_g(G)$.
%\end{myobservation}
%
The following theorem, determining the value of $s_g(G)$ for every connected graph $G$ of order $n\geq 3$, was proved by Anholcer, Cichacz and Milani\v{c} \cite{ref_AnhCic1}.

\begin{mytheorem}[\cite{ref_AnhCic1}]\label{AnhCic1}
Let $G$ be an arbitrary connected graph of order $n\geq 3$. Then
$$
s_g(G)=\begin{cases}
n+2,&\text{if   } G\cong K_{1,3^{2q+1}-2} \text{  for some integer   }q\geq 1,\\
n+1,&\text{if   } n\equiv 2 \imod 4 \wedge G\not\cong K_{1,3^{2q+1}-2} \text{  for any integer   }q\geq 1,\\
n,&\text{otherwise.}
\end{cases}
$$
\end{mytheorem}

%In \cite{ref_Jon} it was proved in turn that for every connected graph $G$ of order $n\geq 3$ 
%$$
%k(G)=\left\{\begin{array}{lll}
%n,&\text{if}& n \not \equiv 2 \imod 4, \\
%n+1,&\text{if}& n\equiv 2 \imod 4.
%\end{array}
%\right.
%$$
%

%[MARCIN] TO WYWALIÅ?EM, BO PRZECIEÅ» JUÅ» Z TWIERDZENIA POWYÅ»EJ WYNIKA, Å»E $n$ NIE ZAWSE WYSTARCZA...

%In order to distinguish $n$ vertices in arbitrary (not necessarily connected) graph we need at least $n$ distinct
%elements of $\gr$. However, $n$ elements are not always enough, as shows the
%following lemma.
%
%\begin{mylemma}[\cite{ref_AnhCic2}]\label{lemma_below}
%Let $G$ be a graph of order $n$. If $n \equiv 2 \imod 4$, then there is no $\gr$-irregular labeling of $G$ for any Abelian group $\gr$ of order $n$.
%\end{mylemma}

Irregularity strength of disconnected graphs was considered in \cite{ref_AnhCic2,ref_AnhCicJurMar}. In particular, 
the following exponential upper bound on $s_g(G)$ was provided there.
%it was proved that $s_g(G)$ is bounded from above.

\begin{mytheorem}[\cite{ref_AnhCicJurMar}]\label{thm_exp}
Let $G$ be a graph of order $n$ having $m$ components, none of which has order less than $3$ and let $p$ be the smallest  number greater than $2^{n-m-1}$ that has all primes distinct in its factorization. Then $s_g(G)\leq p$.
\end{mytheorem}
In the same paper, the authors also presented the exact values and bounds on $s_g(G)$ for disconnected graphs with no star components.

%
%\begin{mytheorem}[\cite{ref_AnhCic2}]\label{modulo}
%Let $G$ be a graph of order $n$ having neither component of order less than $3$ nor a $K_{1,2u+1}$ component for any integer $u\geq 1$. Then:
%
%$$
%\begin{array}{lll}
%k(G)=n,& \text{if} &n \equiv 1\imod 2,\\
%k(G)=n+1,& \text{if} &n\equiv 2 \imod 4,\\
%k(G)\leq n+1, & \text{if} &n\equiv 0 \imod 4.
%\end{array}
%$$
%Moreover for every odd integer $t\geq k(G)$ there exists a $\zet_t$-irregular labeling of $G$.
%\end{mytheorem}
%
%In this paper we give an upper bound for group irregularity strength of all graphs.
%Moreover we give the exact values and bounds on $s_g(G)$ for disconnected graphs with no star components. 

Tutte's zero-flow graph  conjectures, known
as the 5-flow, 4-flow, and 3-flow conjectures, are substantial and seminal sources %a major source
 of inspiration in graph theory \cite{ref_BarTho,ref_JaeLin,ref_Sey,ref_ThoWuZha}.  In these we want to label the edges using non-zero elements of a cyclic group $\zet_k$ so that the sum of the labels flowing into each vertex is equal to the sum of values flowing out of this vertex. In all the above $\gr$-labelings, all the elements of $\gr$ could be used as edge labels as well as the weighted degrees. However, in the case of ordinary irregularity strength, it is forbidden to use $0$ as a label, since admitting it significantly simplifies the problem in many cases, and in fact allows us to focus on an arbitrarily chosen spanning  subgraph of a given graph.  Therefore a natural question arises, namely what is the minimum order of a  group, that allows the \textit{nowhere-zero} (i.e. avoiding the identity element as an edge label) group-irregular labeling. The concept on nowhere-zero modular edge-gracefulness was already considered in \cite{ref_JonZha}. The following open problem was stated in \cite{ref_AnhCic2}. For any Abelian group $\gr$, let $\gr^*=\gr\setminus\{0\}$.
\begin{myproblem}[\cite{ref_AnhCic2}]
Let $G$ be a simple graph with no components of order less than $3$. 
%For any Abelian group $\gr$, let $\gr^*=\gr-\{0\}$. 
Determine the nowhere-zero group irregularity strength ($s_g^*(G)$) of $G$, i.e., the smallest positive integer $k$ such that for any Abelian group $\gr$ of order $k$, there exists a function $f\colon E(G)\rightarrow \gr^*$  such that the sum of edge labels at every vertex is distinct.
\end{myproblem}
Obviously $s_g(G)\leq s_g^*(G)$ and this equality is sharp,  as  the example of $K_{1,4n-1}$ shows \cite{ref_JonZha}.  In this paper we give a linear upper bound on these graph invariants, showing these are not greater than twice the order of given graph, see Corollary~\ref{GroupIrregCor} below.

%%%%%%%%%%%%%%%%%%

Not less intriguing problem than the concept of irregularity strength itself is its local version, where we require only distinction between weights at adjacent vertices. It is believed that just labels $1,2,3$ are sufficient to achieve this for every graph without a component of order $2$, see~\cite{Louigi30,Louigi,KalKarPf_123,123KLT,ref_ThoWuZha,123with13} for results related to this well-known and widely studied so-called \textit{1--2--3 Conjecture}, which remains open in general.  
%$1,2,3$ See... for results related with this concept, and in particular..., where it is proved that weights $1,2,3,4,5$ %are sufficient for every graph without component of order $2$  results of graphs 1--2--3 Conjecture
Analogously, for any given Abelian group $\gr$, by a vertex-coloring $\gr$-labeling of a given graph $G$ we understand a labeling of its edges with elements of $\gr$ which induces a proper vertex coloring with the resulting sums at vertices.
%One can also focus on distinguishing the weighted degrees only locally or, in other words, on choosing $\gr$-labeling that induces a proper coloring of $G$ with colors corresponding with weighted degrees, which is called a \textit{vertex-$\gr$-coloring}. 
The smallest $k$ admitting such a labeling for any group $\gr$ of order $k$ %order of a group that allows such a labeling 
is called the \textit{group sum chromatic number} and denoted by $\gchi(G)$. Anholcer and Cichacz in \cite{ref_AnhCic3} proved that for arbitrary graph not having components of order less than $3$, $\gchi(G)\in \{\chi(G),\chi(G)+1\}$ and %they 
completely characterized the graphs for which $\gchi(G)=\chi(G)+1$. However, also in this case there was no restriction imposed on a possible usage of $0$ as an edge label or as a weighted degree, and the following open problem was formulated there.

%[MARCIN]Tu był kompletny burdel w notacji
\begin{myproblem}[\cite{ref_AnhCic3}]
Let $G$ be a simple graph with no components of order less than $3$. 
%For any Abelian group $\gr$, let $\gr^*=\gr-\{0\}$. 
Determine the nowhere-zero group sum chromatic number (${\gchi}^\star(G)$) of $G$, i.e., the smallest positive integer $k$ such that for any Abelian group $\gr$ of order $k$, there exists a function $f\colon E(G)\rightarrow \gr^*$  such that the resulting sums at the vertices properly color them. %of edge labels properly color the vertices.
\end{myproblem}
Obviously $\gchi(G)\leq {\gchi}^\star(G)$ and this equality is sharp. For instance, $\gchi(C_{4n})=2$, as showed in \cite{ref_AnhCic3}, while one can easily see that ${\gchi}^\star(C_{4n})>2$. 
%Lub G=C_4n

In this article we provide an upper bound of twice the %double 
maximum degree for ${\gchi}^\star(G)$, see Theorem~\ref{LocalStarUpperBoundLem} for more details.
Moreover we prove a finite upper bound for trees, see Theorem~\ref{tree2}, and more generally -- for any family for graphs with bounded arboricity, in particular for planar graphs, see Subsection~\ref{ArboricitySubsection}.   
% give an upper bound on the order of a group that guarantees the existence of nowhere-zero vertex-$\gr$-coloring  %of $G$. 

%%%%%%%%%%%%%%%%%%%%%%%%%%%%%%%%%%

%[MARCIN]Przeniosłem wnioski tutaj i dopisałem krótkie wstępy. Generalnie uważam, że twierdzenia i wnioski powinny pójść do wstępu, a dowody jako osobne sekcje.

\section{Nowhere-zero group irregularity strength}

We first prove an upper bound for the nowhere-zero group irregularity strength of a graph $G$ linear in terms of the order of $G$. We will moreover prove this result under the condition that the identity element of a given group $\gr$ cannot be induced as the weighted degree of any vertex of $G$ -- we call such a $\gr$-labeling $f$ of $G$ \textit{non-zero} (as the induced weighting function $w_f$ cannot take value $0$).

%First of the main results of this article is the following theorem:

\begin{mytheorem}\label{thm_main}
Let $G$ be arbitrary graph of order $n$ having no component of order less than $3$. Let $\gr$ be arbitrary Abelian group of order at least $2n$. Then there exists a $\gr$-irregular labeling of $G$, such that no edge is labeled with $0$ and no vertex has vertex degree equal to $0$.
\end{mytheorem}

\textit{Proof.} The proof follows by induction in the number of edges.

Suppose first that $G$ is a path %Let us start with smallest possible graph $G$, i.e. 
$P_{3}$ with vertices, say, $u$, $v$ and $w$ and edges $uv$ and $uw$. Let $\gr$ be arbitrary Abelian group of order at least $2|G|=6$. Choose any element $a\in \gr\setminus\{0\}$ and set $f(uv)=a$. Now, choose any $b\in\gr\setminus\{0,a,-a\}$ and set $f(uw)=b$. Both edge labels are different than $0$, and so are the vertex weighted degrees, since $w(u)=a+b$, $w(v)=a$ and $w(w)=b$. It is also obvious that the weighted degrees are three distinct elements of $\gr$. Note that the choice of such $a$ and $b$ is always possible if $\gr$ has at least $4$ elements, and it is the case, since $|\gr|\geq 6$.

Now let $G$ be arbitrary graph of order $n$ with at least $3$ edges having no component of order less than $3$ and let $\gr$ be any Abelian group of order at least $2n$. In the induction step we can assume that for every proper subgraph $H$ of $G$ having no component of order less than $3$ and for every Abelian group $\gr'$ of order at least $2|H|$, there is a $\gr'$-irregular labeling $f_H$ of $H$ in which no edge has label $0$ and no vertex has weighted degree $0$. In particular, there is such labeling of $H$ with $\gr^\prime=\gr$, since $|\gr|\geq 2n\geq 2|H|$. We will extend $f_H$ to the labeling $f$ of $G$, having the same properties.

We choose $H$ in one of the following ways. If there is a component $C\cong P_3$ of $G$, then $H=G-C$. Otherwise, if there is a component $C$ and an edge $e\in E(C)$ not being a bridge in $C$, then $H=G-e$. Finally, if $G$ is a forest with each component of order at least $4$, then choose any leaf edge $e$ of any component and let $H=G-e$.

Let us consider the first case. Assume that $G=H\cup H^\prime$, where $V(H^\prime)=\{u,v,w\}$ and $E(H^\prime)=\{uv,uw\}$. Let $f_H$ be a nowhere-zero and non-zero $\gr$-irregular labeling of $H$, existing by the induction hypothesis. Now let $f(e)=f_H(e)$ for $e\in E(H)$. Now choose any element of $a\in \gr$ such that $a\neq 0$ and $a\neq w(x)$ for $x\in V(H)$ and set $f(uv)=a$. Such $a$ can be chosen, as only $n-3$ vertex weighted degrees have been assigned so far and $|\gr|>n-2$. Now choose $b\in\gr$ such that $b\not\in \{0,a,-a\}$, $b\neq w(x)$ for $x\in V(H)$ and $b\neq w(x)-a$ for $x\in V(H)$ and set $f(uw)=b$. The number of forbidden elements is equal to at most $3+2(n-3)=2n-3<|\gr|$, so we can choose such $b$. Obviously, the two new edge labels are not $0$ and neither are the three new weighted degrees $w(u)=a+b$, $w(v)=a$ and $w(w)=b$. Also, the three new weighted degrees are pairwise distinct and not equal to any $w(x)$, where $x\in V(H)$. Thus $f$ is a nowhere-zero and non-zero $\gr$-irregular labeling of $G$.

In the second case, let $H=G-e$ and let $f_H$ be a nowhere-zero and non-zero $\gr$-irregular labeling of $H$, guaranteed by the induction hypothesis. Now let $f(y)=f_H(y)$ for $y\in E(H)$. Let us denote the vertices incident with $e$ in $G$ by $u$ and $v$. Choose an element $a\in\gr$ such that $a\not\in \{0,-w_H(u),-w_H(v)\}$, $a\neq w_H(x)-w_H(u)$ for $x\in V(G)\setminus \{u,v\}$ and $a\neq w_H(x)-w_H(v)$ for $x\in V(G)\setminus \{u,v\}$. Set $f(e)=a$. The number of forbidden values is at most $3+2(n-2)=2n-1<|\gr|$, so we can always choose such $a$. Note that two adjusted weighted degrees remain distinct and because of the way that $a$ was chosen, they are different than any weighted degree $w(x)$ for $x\in V(G)\setminus \{u,v\}$. This means that $f$ is a nowhere-zero and non-zero $\gr$-irregular labeling of $G$.

Finally, consider the third case. Assume that the ends of $e$ are $u$ and $v$, where $u$ is the pendant vertex. Having a nowhere-zero and non-zero $\gr$-irregular labeling $f_H$ of $H$, we set $f(y)=f_H(y)$ for $y\in E(H)$. Then we choose $a\in \gr$ such that $a\not\in \{0,-w_H(v)\}$,  $a\neq w_H(x)$ for $x\in V(G)\setminus \{u,v\}$ and $a\neq w_H(x)-w_H(v)$ for $x\in V(G)\setminus \{u,v\}$. There are at most $2+2(n-2)=2n-2<|\gr|$ forbidden values, so we can choose such $a$. The adjusted weighted degree $w(v)$ and the new weighted degree $w(u)$ are distinct and different than any of the weighted degrees $w(x)$ for $x\in V(G)\setminus \{u,v\}$, so also in this case $f$ is a nowhere-zero and non-zero $\gr$-irregular labeling of $G$. This completes the proof.
\qed

\ \\
%Note that we can obtain slightly better bounds if we do not restrict ourselves to nowhere-zero flows. To be more specific, the necessary size of the group in the base case of the induction would be $3$ (instead of $4$), and in the induction step, the number of forbidden elements would be respectively $2n-5$, $2n-4$ and $2n-4$. 
The following corollary immediately follows.

\begin{mycorollary}\label{GroupIrregCor}
Let $G$ be arbitrary graph of order $n$ having no component of order less than $3$. Then
$$
s_g(G)\leq s_g^*(G)\leq 2n.
$$
\end{mycorollary} 

Note that this significantly improves the result in Theorem \ref{thm_exp}.

\section{Nowhere-zero group sum chromatic number} %of graphs}

We will make use of the notion of the \textit{coloring number} of a graph, introduced by  Erd\H{o}s and Hajnal in \cite{ErdHaj}. For a given graph $G$ by ${\rm col}(G)$ we denote its coloring number, that is the least integer $k$ such that each subgraph of $G$ has minimum degree less than $k$.
Equivalently, it is the smallest $k$ for which we may linearly order all vertices of $G$ into a sequence $v_1,v_2,\ldots,v_n$ so that every vertex $v_i$ has at most $k-1$ neighbors preceding it in the sequence.
Hence $\chi(G)\leq {\rm col}(G)\leq \Delta(G)+1$.
Note that ${\rm col}(G)$ equals the degeneracy of $G$ plus $1$, and thus the result below may %also
be formulated %by means of such a well-known graph invariant instead of ${\rm col}(G)$.
in terms of either of the two graph invariants.

\subsection{General Upper Bound}

In order to prove a general upper bound for the nowhere-zero group sum chromatic number we will need to consider the following more general setting.

%In the general case if we omit the assumption that all vertices must have  nonzero weighted degree we obtain the %following theorem.
%
% [MARCIN]W lemacie chyba nie trzeba zakładać spójności?
%Poza tym, zwiększyłbym rangę na twierdzenie.
%To zmienia również wnioski (po co spójność?)
%Na poczÄtku dowodu dodałem zdanie o tym, że można się skoncentrować na grafach spójnych.

\begin{mytheorem}\label{LocalStarUpperBoundLem}
For every graph $G=(V,E)$ of order at least $3$, a set of marked vertices $M\subset V$ and an abelian group $\gr$ with $|\gr|\geq \Delta(G)+{\rm col}(G)-1$, there exists a labeling  $f\colon E(G)\rightarrow \gr^*$ so that $w(v)\neq 0$ for every $v\in M$ and $w(u)\neq w(v)$ for every edge $uv\in E$ such that $u,v\notin M$.
\end{mytheorem}

{\it Proof.} Since we consider local distinguishing of vertices, it is enough to prove the result for connected graphs. It is straightforward to verify it in the case when $G$ has $3$ vertices. So we assume that $|V|\geq 4$ and prove the theorem by induction with respect to $|V|$. 

Let $\gr$ be any abelian group with $|\gr|\geq \Delta(G)+{\rm col}(G)-1$, and let $M$ be any set of marked vertices in $G$. Suppose $v$ is a vertex of minimum degree $\delta$ in $G$, and let $N_G(v)=\{v_1,v_2,\ldots,v_\delta\}$ (note that $\delta\leq {\rm col}(G)-1$). %Set $G'=G-v$.

If $\delta=1$, set $G'=G-v$. Then $G'$ is connected, has order at least $3$, and $\Delta(G')+{\rm col}(G')-1 \leq \Delta(G)+{\rm col}(G)-1$, so we may label $G'$ with the elements of $\gr^*$ by induction (consistently with the thesis of the theorem) with the set of marked vertices $M':=M\cup\{v_1\}\setminus\{v\}$. In order to extend the obtained labeling $f_{G'}$ of $G'$ to a labeling $f$ of the whole $G$, we set $f(e)=f_{G'}(e)$ for every $e\in E\setminus\{v_1v\}$, and choose a label $f(v_1v)\in \gr$ so that $f(v_1v)\neq 0$ and so that $v_1$ is sum distinguished from its (at most $\Delta(G)-1\geq 1$) neighbors other than $v$ if $v_1\notin M$ (since $v_1\in M'$, it is guaranteed that $w(v)\neq w(v_1)$), or so that $w(v_1)\neq 0$ if $v_1\in M$. We can do it, as $|\gr|> \Delta(G)$. It is easy to verify that the obtained $f$ fulfills our requirements (regardless of the fact whether $v$ is in $M$ or not). 

We may thus assume that $\delta\geq 2$ (hence ${\rm col}(G)\geq 3$). %Then let $I$ be the set of
Then set $G'=G-\{v_1v,v_2v\}$, and let $G_1,\ldots,G_k$ be the components of $G'$ (hence $k\leq 3$).
Since $v$ has minimum degree in $G$ and $|V|\geq 4$, no component of $G'$ has order $2$. %at least $3$ 
Therefore, similarly as above, by induction, for each $i\in\{1,\ldots,k\}$ there exists a labeling $f_{G_i}:E(G_i)\rightarrow \gr^*$ of $G_i$, consistent with the thesis of the theorem, with the set of marked vertices $M_i=M\cap V(G_i)\setminus\{v\}$.

We define $f:E\rightarrow \gr^*$ by first setting $f(e)=f_{G_i}(e)$ for every $e\in E(G_i)$, $i=1,\ldots,k$. Then we choose a label $f(v_2v)\in \gr$ so that $f(v_2v)\neq 0$ and $w(v)\neq w(v_1)$ (which cannot be influenced by a later choice of $f(vv_1)$, as this counts in the sums of both, $v$ and $v_1$), and so that %$w(v_2)\neq 0$ and 
$v_2$ is sum distinguished from its (at most $\Delta(G)-1\geq 1$) neighbors other than $v$
(and $v_1$ if $v_1v_2\in E$) if $v_2\notin M$, or so that $w(v_2)\neq 0$ if $v_2\in M$. 
We can do it, as $|\gr|\geq \Delta(G)+{\rm col}(G)-1> 1+1+(\Delta(G)-1)$.
Finally we choose a label $f(v_1v)\in \gr$ so that $f(v_1v)\neq 0$ and so that
$v_1$ is sum distinguished from its (at most $\Delta(G)-1\geq 1$) neighbors other than $v$
if $v_1\notin M$, or so that $w(v_1)\neq 0$ if $v_1\in M$, and moreover so that
$v$ is sum distinguished from its (at most $\delta-1\geq{\rm col}(G)-2$) neighbors other than $v_1$
if $v\notin M$, or so that $w(v)\neq 0$ if $v\in M$. 
We can do it, as $|\gr| > 1+ (\Delta(G)-1)+({\rm col}(G)-2)$.
It straightforward to verify that the obtained $f$ complies with our requirements.~\qed\\

Setting $M=\emptyset$, we immediately obtain the following result.

%[MARCIN]Problemy z notacją powtórzyły się tutaj
\begin{mycorollary}\label{LocalStarUpperBoundCor}
If $G$ is a graph with no component of order less than $3$, then ${\gchi}^\star(G)\leq \Delta(G)+{\rm col}(G)-1\leq 2\Delta(G)$.
\end{mycorollary}

Taking into account that for every planar graph $G$ we have ${\rm col}(G)\leq6$, we thus immediately obtain for instance the following corollary.

%[MARCIN]I tutaj
\begin{mycorollary}\label{LocalStarUpperBoundCorPlanar}
If $G$ is a planar graph with no component of order less than $3$, then ${\gchi}^\star(G)\leq \Delta(G)+5$.
\end{mycorollary}

Note also that if we additionally want to forbid zero sums at all vertices, then within the proof of Theorem~\ref{LocalStarUpperBoundLem} above, we obtain at most two additional constraints while choosing a label for a given edge (i.e. forbidden zero-sums at its ends). Consequently, by %the analogon 
a straightforward adaptation 
of the proof above, we obtain the following. 

\begin{myobservation}\label{chromatic}
Let $G$ be arbitrary graph of order $n$ having no component of order less than $3$. Let $\gr$ be arbitrary Abelian group of order at least $\Delta(G)+{\rm col}(G)+1$. Then there exists a vertex--coloring $\gr$-labeling of $G$ such that no edge is labeled with $0$ and no vertex has weighted degree equal to $0$.
\end{myobservation}

\subsection{Trees}

\begin{mytheorem}\label{tree2}
%EXTEND:
%If $T$ is a tree of order $n\geq 3$, then ${\gchi}^\star(G)\leq 4$.
Let $T$ be arbitrary tree of order $n\geq 3$. Let $\gr$ be arbitrary Abelian group of order at least $4$. Then there exists a vertex-coloring $\gr$-labeling of $T$ such that no edge is labeled with $0$.
\end{mytheorem}

{\it Proof.} 
For $n=3$ the thesis obviously holds, so we may assume that $n\geq 4$ and prove the theorem by induction with respect to $n$. 

Let $\gr$ be an Abelian group with $|\gr|\geq 4$.
Root $T$ at some leaf $r$, and let $P$ be a maximal path in $T$ starting at $r$. Let $v$ be the second last (counting from $r$) vertex on this path; denote its degree by $d$. Then all sons of $v$ are leaves -- denote them by $v_1,v_2,\ldots,v_{d-1}$. Let $u$ be the father of $v$ in $T$.

Let $T'=T-\{v_1,v_2,\ldots,v_{d-1}\}$. By induction there exists  
a $\gr^*$-labeling $f_{T'}$ %of $T'$ 
inducing distinct sums for the neighbours in $T'$, unless $T'$ is an isolated edge -- we color such edge with any non-zero label from $\gr$ then. We will extend $f_{T'}$ to a desired labeling $f$ of the entire $T$.

Suppose first that there exists $c\in\gr^*$ such that $(d-2)c=0$. Then we assign $c$ to $vv_1,vv_2,\ldots,vv_{d-2}$ (if there are any). Note that the obtained temporary sum at $v$ is non-zero, hence in order to finish a labeling of $T$ it is sufficient to choose a non-zero label $f(vv_{d-1})$ so that $w(v)\neq c$ and $w(v)\neq w(u)$ afterwards -- this is always feasible, as $|\gr|\geq 4$.

On the other hand, if we have $(d-2)c'\neq 0$ for every $c'\in\gr^*$ (thus $d\geq 3$), then it is straightforward to notice that $(d-2)c'\neq(d-2)c''$ for every two distinct elements $c',c''\in\gr$. Therefore, there exists $c\in\gr^*$ such that the temporary sum at $v$ will be non-zero after setting $f(vv_1)=c, f(vv_2)=c,\ldots,f(vv_{d-2})=c$. Thus we may finalize the  labeling of $T$ by the same argument as above (i.e., by choosing non-zero $f(vv_{d-1})$ so that $w(v)\neq c$ and $w(v)\neq w(u)$). 
~\qed

\begin{mycorollary}\label{tree3}
If $T$ is a tree of order $n\geq 3$, then ${\gchi}^\star(G)\leq 4$.
\end{mycorollary}

Note that the obtained bound for trees is tight, as not for every tree there exists a $\gr^*$-labeling inducing distinct sums for the neighbours for $\gr\cong \zet_3$. To see this consider $e.g.$ a symmetric double star with $5$ edges, i.e. a graph of maximum degree $3$ which might be obtained from two separate stars $K_{1,3}$ by identifying one edge from the first one with any edge from the second one. 
%(note that in order to distinguish neighbours of degree $3$ in such a tree, we must put $1$ and $2$ on the two %pendant edges edges incident with the sa...NIEPRAWDA, potrzeba nieco dluzszego uzasadnienia). 

Analogously as in the previous subsection, by the same reasoning as in the proof of Theorem~\ref{tree2} above, we may obtain the following result %conclude the following 
if we additionally forbid zero sums at all vertices %assume that 
(this time, while choosing a label for a given edge, we have merely at most one additional constraint -- so that $v$ receives a non-zero sum).

\begin{myobservation}\label{tree}
Let $T$ be arbitrary tree of order $n\geq 3$. Let $\gr$ be arbitrary Abelian group of order at least $5$. Then there exists a vertex-coloring $\gr$-labeling of $T$, such that no edge is labeled with $0$ and no vertex has weighted degree equal to $0$.
\end{myobservation}

This is also tight.
Indeed, observe that for $\gr\cong \zet_3$ or  $\gr\cong\zet_2\times \zet_2$ there does not exists 
a vertex-coloring $\gr$-labeling of $K_{1,3}$, such that no edge is labeled with $0$ and no vertex has weighted degree equal to $0$.

%, therefore the above bound for trees is tight. However the bound can be slightly improved if we make weaker assumptions and allow vertices to have weighted degree equal to $0$.\\

In the next subsection we in particular eventually show a finite upper bound for ${\gchi}^\star(G)$ for the family of planar graphs, 
by providing a %more 
general upper bound for this graph invariant in terms of so-called arboricity of a graph.
%
%by more generally providing ...zwiazujac... z so-called arboricity of graphs. 
%The results of this subsection will be key ingredients of our reasoning in the next one.
One of the key ingredients of the proof of this fact will be the results from this subsection on trees.

\subsection{Graphs with Bounded Arboricity, Planar Graphs}\label{ArboricitySubsection}

By the \textit{arboricity}, $a(G)$ of a graph $G$ we mean the least number of forests into which we may decompose the set of edges of $G$, i.e., in other words, the least number of colors with which we can color the edges of $G$ so that each color induces a forest in $G$, see~\cite{Nash-Wiliams}. 

\begin{mylemma}\label{ForestDecompLemma}
Suppose $G$ is a graph without isolated edges and isolated triangles. Then we may color the edges of $G$ with $a(G)$ colors so that each color induces a forest without isolated edges.
\end{mylemma}

{\it Proof.} 
We may assume that $G$ is connected and has at least $4$ edges.
Consider a coloring of the edges of $G$ with $a(G)$ colors with the least number of monochromatic components being isolated edges, and suppose this number is positive. Let $uv$ be an isolated edge colored, say, with color $1$ (isolated means isolated in the forest colored with the same color).

We first observe that every other monochromatic component incident with $u$ or $v$ must contain exactly $2$ edges -- one incident with $u$ and the other incident with $v$.
In order to see this consider a monochromatic component $T$ incident with $u$ or $v$ 
%which has $1$ or at least $3$ edges and is 
colored, say, with $2$.
Note that $T$ must contain a path joining $u$ with $v$, as otherwise we could recolor $uv$ using color $2$ (note that no monochromatic cycle would be created in such a case), and thus reduce the number of  isolated edges in $G$, a contradiction. So suppose now that $T$ has $1$ or at least $3$ edges. 
Then however, $T$ cannot contain a pendant edge $e$ adjacent with $uv$, as otherwise we could recolor $e$ coloring it with $1$, and this way reduce the number of isolated edges in $G$. Consequently, 
 $T$ must have at least $6$ edges, and thus we may recolor with $1$  one of the edges of $T$ adjacent with $uv$ which lays on the path joining $u$ with $v$ in $T$, and this way obtain a contradiction with minimality of the number of isolated edges in monochromatic components of $G$.

%, for suppose first that $uv$ is adjacent with some other monochromatic isolated edge colored, say, with $2$. Then %we can recolor $uv$ coloring it with $2$, and thus reducing the number of monochromatic isolated edges, a %contradiction. Suppose then that $uv$ is incident with a monochromatic component, say $T$, with at least $3$ %edges colored, say, with $2$. Then we may similarly as above recolor $uv$ using color $2$, unless this creates a %monochromatic cycle (colored $2$). Thus suppose there is a path $P$ joining $u$ and $v$ in $T$. If $uv$ is %adjacent with a pendant edge of $T$, say $e$, then we may recolor $e$ using $1$ in order to reduce the number of %monochromatic isolated edges in $G$, and thus obtaining a contradiction. Therefore, $T$ must have at least $6$ %edges, and we may recolor one of the edges of $T$ adjacent with $uv$ using $1$, and thus obtaining a %contradiction.

Consequently, as $G$ is connected and has at least $4$ edges, $uv$ must be incident with two monochromatic components of size two, each containing both $u$ and $v$, colored differently, say a path $P_1=uxv$ colored with $2$  and a path $uyv$ colored with $3$.  Then however we may recolor $vx$ using $3$ and $uv$ using $2$, reducing the number of  isolated edges, and thus obtaining a contradiction. 
~\qed\\

\begin{mycorollary}\label{GroupArboricity1}
Let $G$ be a graph of arboricity at most $a$ containing no isolated edges. Let $\gr\cong\gr_1\times \gr_2\times \ldots \times \gr_a$ be an Abelian group being the product of Abelian groups $\gr_i$ such that $|\gr_i|\geq 4$ for $i=1,2,\ldots,a$. Then there exists a vertex-coloring $\gr$-labeling of $G$ such that no edge is labeled with $0$.
\end{mycorollary}

{\it Proof.} 
We may assume that $G$ is connected. The theorem obviously holds if $G$ is a triangle, so we may assume this is not the case. Consequently, by Lemma~\ref{ForestDecompLemma}, the edges of $G$ can be colored with colors $1,2,\ldots,a$ so that each color induces a forest without isolated edges -- denote these forests by $F_1,F_2,\ldots,F_a$, respectively. By Theorem~\ref{tree2}, for every $i\in \{1,2,\ldots,a\}$ there exists a nowhere-zero vertex-coloring $\gr_i$-labeling $f_i$ of $F_i$, i.e. $f_i:E(F_i)\to \gr_i$. We extend each of these to an edge labeling of $G$: 
$g_i:E(G)\to\gr_i$ by setting $g_i(e):=f_i(e)$ if $e\in F_i$ or $g_i(e)=0$ otherwise, for $i=1,2,\ldots,a$. Now let us define $f:E(G)\to \gr_1\times \gr_2\times \ldots \times \gr_a$ by setting $f(e)=(g_1(e),g_2(e),\ldots,g_a(e))$ for every $e\in E(G)$. To see that $f$ is a  a nowhere-zero vertex-coloring $\gr$-labeling of $G$, consider any edge $uv\in E(G)$. Obviously $uv\in E(F_j)$ for some $j\in\{1,2,\ldots,a\}$. Then $w_f(u)=(w_{g_1}(u),w_{g_2}(u),\ldots,w_{g_a}(u))$ and $w_f(v)=(w_{g_1}(v),w_{g_2}(v),\ldots,w_{g_a}(v))$, and hence $w(u)\neq w(v)$, as by the definiton: $w_{g_j}(u)=w_{f_j}(u)\neq w_{f_j}(v)=w_{g_j}(v)$.  
~\qed\\

\begin{mycorollary}\label{GroupArboricity2}
Let $G$ be a graph of arboricity at most $a$ containing no isolated edges and let $k(a)$ be the least positive integer such that each Abelian group $\gr$ of order $k_a$ is isomorphic to some product $\gr_1\times \gr_2\times \ldots \times \gr_a$ of Abelian groups $\gr_i$ such that $|\gr_i|\geq 4$ for $i=1,2,\ldots,a$. Then ${\gchi}^\star(G)\leq k(a)$.
\end{mycorollary}

The fundamental theorem of finite Abelian groups states that a finite Abelian group $\gr$ of order $n$ can be expressed as the direct product of cyclic subgroups of prime-power orders. This implies that
$$\gr\cong\zet_{p_1^{\alpha_1}}\times\zet_{p_2^{\alpha_2}}\times\ldots\times\zet_{p_k^{\alpha_k}}\;\;\; \mathrm{where}\;\;\; n = p_1^{\alpha_1}\cdot p_2^{\alpha_2}\cdot\ldots\cdot p_k^{\alpha_k}$$
and $p_1,p_2,\ldots,p_k$ are not necessarily distinct primes, see e.g.~\cite{ref_Gal}. Recall also that for any positive integer $i$ there exists a prime number between $i$ and $2i$, see \cite{ref_The}.
Therefore, for any fixed $a$, the value of $k(a)$ from Corollary~\ref{GroupArboricity2} exists (is finite) and is e.g. upper bounded by $p_1\cdot p_2\cdot\ldots\cdot p_a$, where $p_1,p_2,\ldots,p_a$ are first consecutive (pairwise different) $a$ prime numbers larger than $3$.

\begin{mycorollary}\label{GroupArboricity3}
Let $G$ be a graph of arboricity at most $a$ containing no isolated edges and let $p_1,p_2,\ldots,p_a$ be the first consecutive  $a$ prime numbers larger than $3$. Then ${\gchi}^\star(G)\leq p_1\cdot p_2\cdot\ldots\cdot p_a$.
\end{mycorollary}

The upper bound from Corollary~\ref{GroupArboricity3} above can be slightly improved via direct application of Corollary~\ref{GroupArboricity2} and careful analysis of powers of $2$ and $3$. We exemplify this in the setting of planar graphs within the following proof. 
%
%setting
%
%analysis of the case of planar graphs. 
Recall the well-known result of Nash-Williams, \cite{Nash-Wiliams}, for a graph $G$:
\begin{equation}\label{Nash-Williams_bound}
a(G) = \max_{H\subseteq G} \left\lceil\frac{|E(H)|}{|V(H)-1|}\right\rceil.
\end{equation}

\begin{mycorollary}\label{GroupArboricityPlanar}
For every planar graph $G$ containing no isolated edges, ${\gchi}^\star(G)\leq 140$.
\end{mycorollary}

{\it Proof.} As due to Euler's Formula, every planar graph with $n$ vertices has at most $3n-6$ edges, by~(\ref{Nash-Williams_bound}), $a(G)\leq 3$. Let $\gr$ be an Abelian group with $|\gr|=140$. Note that $140=2^2\cdot5\cdot 7$. Thus one of the following must hold:
$$\gr\cong\zet_{2^{2}}\times\zet_{5}\times\zet_{7}\;\;\; {\rm or} \;\;\; \gr\cong\left(\zet_{2}\times\zet_2\right)\times\zet_{5}\times\zet_{7}.$$
By Corollary~\ref{GroupArboricity2} we thus obtain the thesis.
~\qed\\

%\begin{myconjecture}There exists a constant $C$ such that if $G$ is a planar graph with no component of order %less than $3$, then  ${\gchi}^\star(G)\leq C$.
%\end{myconjecture}

\section{Final remarks}

We conclude the paper by posing the following conjectures.

\begin{myconjecture}
There exists a constant $C$ such that
for every graph $G$ of order $n$ having no component of order less than $3$:
%Let $G$ be arbitrary graph of order $n$ having no component of order less than $3$. Then there exists a constant %$C$ such that
$$
s_g^*(G)\leq n+C.
$$
\end{myconjecture}

\begin{myconjecture}
There exists a constant $C$ such that
for every graph $G$ having no component of order less than $3$:
%Let $G$ be arbitrary graph of order $n$ having no component of order less than $3$. Then there exists a constant %$C$ such that
$${\gchi}^\star(G)\leq \chi(G) +C.$$
\end{myconjecture}

In particular we state the following for planar graphs.

\begin{myconjecture}If $G$ is a planar graph with no component of order less than $3$, then  $${\gchi}^\star(G)\leq 8.$$
\end{myconjecture}

\nocite{*}
\bibliographystyle{amsplain}

\end{document}